\documentclass[preprint,12pt]{elsarticle}
\usepackage{geometry}
\usepackage{listings}
\usepackage{color}

\usepackage{tikz}
\usetikzlibrary{backgrounds}
\usepackage{graphicx}
\usepackage{bm}
\usepackage{amsmath}
\usepackage{setspace}
\usepackage{hyperref}
\usepackage{subfigure}
\usepackage{multirow}
\hypersetup{
    colorlinks=true,       
}
\usepackage{amssymb}

\usepackage{lineno}
\usepackage{diagbox}

\usepackage{amsmath}
\usepackage{ntheorem}

\usepackage{algorithm}
\usepackage{algorithmic}
\usepackage{changes}
\usepackage{soul}  

\def \ux {\frac{\partial u}{\partial x}}
\def \ut {\frac{\partial u}{\partial t}}
\def \uxx {\frac{\partial^2 u}{\partial x^2}}

\begin{document}
\begin{frontmatter}


\title{Sparse Deep Neural Network for Nonlinear Partial Differential Equations}



\author[label1]{Yuesheng Xu}
\ead{y1xu@odu.edu}
\author[label2,label3]{Taishan Zeng\corref{cor1}}
\address[label1]{Department of Mathematics and Statistics, Old Dominion University,\\ Norfolk, Virginia 23529, USA}
\address[label2]{School of Mathematical Science, South China Normal University,\\ Guangzhou 510631, China}
\address[label3]{Guangdong Key Laboratory of Big Data Analysis and Processing,\\  Guangzhou 510006, China}
\ead{zengtsh@m.scnu.edu.cn}
\cortext[cor1]{Corresponding author}
\fntext[label4]{
Y. Xu is supported in part by US National Science Foundation under grant DMS-1912958. 
T. Zeng is supported in part by the National Natural Science Foundation of China under grant 12071160 and U1811464, by the Natural Science Foundation of Guangdong Province under grant 2018A0303130067,  by the Opening Project of Guangdong Province Key Laboratory of Computational Science at the Sun 
Yat-sen University under grant 2021022, and by the Opening Project of Guangdong Key Laboratory of Big Data Analysis and Processing under grant 202101. }

\begin{abstract}
More competent learning models are demanded for data processing due to increasingly greater amounts of data available in applications. Data that we encounter often have certain embedded sparsity structures. That is, if they are represented in an appropriate basis, their energies can concentrate on a small number of basis functions. This paper is devoted to a numerical study of adaptive approximation of solutions of nonlinear partial differential equations whose solutions may have singularities, by deep neural networks (DNNs) with a sparse regularization with multiple parameters. Noting that DNNs have an intrinsic multi-scale structure which is favorable for adaptive representation of functions, by employing a penalty with multiple parameters, we develop DNNs with a multi-scale sparse regularization (SDNN) for effectively representing functions having certain singularities. We then apply the proposed SDNN to numerical solutions of the Burgers equation and the Schr\"odinger equation. Numerical examples confirm that solutions generated by the proposed SDNN are sparse and accurate.

\end{abstract}


\end{frontmatter}



\section {Introduction}
The goal of this paper is to develop a sparse regularization deep neural network model for numerical solutions of nonlinear partial differential equations whose solutions may have singularities.  We will mainly focus on designing a sparse regularization model by employing multiple parameters to balance sparsity of different layers and the overall accuracy.  The proposed ideas are tested in this paper numerically to confirm our intuition and more in-depth theoretical studies will be followed in a future paper. 

Artificial intelligence especially deep neural networks (DNN) has received great attention in many research fields.
From the approximation theory point of view, a neural network is built by functional composition to approximate a continuous function with arbitrary accuracy. Deep neural networks are proven to have better approximation by 
practice and theory due to their relatively large number of hidden layers. Deep neural network has achieved state-of-the-art performance in a wide range of applications, including speech recognition \cite{Dahl2012}, 
computer vision \cite{Krizhevsky_2012}, natural language processing \cite{Bert2018},  and finance \cite{Chong2017}. For an overview of deep learning the readers are referred to monograph \cite{Goodfellow2016}.
Recently, there was great interest in applying deep neural networks to the field of scientific computing, 
such as discovering the differential equations from observed data \cite{Raissi2018}, solving the partial 
differential equation (PDE) \cite{Han_jentzen_E_2018, Lagaris1998, Lagaris2000, PINN_2019}, and problem aroused in physics \cite{Diefenthaler2021}.   
Mathematical understanding of deep neural networks received much attention in the applied mathematics community. A universal approximation theory of neural network for Borel measurable function on compact domain is established in \cite{Cybenko1989}. Some recent research studies the expressivity of deep neural networks for different function spaces \cite{DeVore2021}, for example,  Sobolev spaces, Barron functions, and H{\" o}lder spaces.  There are close connections between deep neural network and traditional approximation methods, such as splines \cite{Daubechies2019, Unser2019}, compressed sensing \cite{Adcock2021}, and finite elements \cite{He_Li_Xu_2020, Jung2020}.  Convergence of deep neural networks and deep convolutional neural networks are studied in \cite{XuZhang2021} and \cite{XuZhang_cnn2021} respectively.    
Some work aims at understanding the training process of DNN. For instance, in paper \cite{Cyr2020}, the training process of DNN is interpreted as learning adaptive basis from data.

Traditionally, deep neural networks are dense and over-parameterized. A dense network model requires more memory and other computational resources during training and inference of the model. 
Increasingly greater amounts of data and related model sizes demand the availability of more competent learning models. Compared to dense models, sparse deep neural networks require less memory, less computing time and have better interpretability. Hence, sparse deep neural network models are desirable. On the other hand, animal brains are found to have hierarchical and sparse structures \cite{Friston2008}. The connectivity of an animal brain becomes sparser as the size of the brain grows larger. Therefore, it is not only necessary but also natural to design sparse networks. In fact, it was pointed out in \cite{Hoefler_2021SIAMNews} that the future of deep learning relies on sparsity. Furthermore, over-parameterized and dense models tend to lead to overfitting and weakening the ability to generalize over unseen examples. Sparse models can improve accuracy of approximation. Sparse regularization is a popular way to learn the sparse solutions \cite{Candes2008, Xu_L0_2021, Xu_Ye2019, Zhang-Xu-Zhang}. The readers are referred to \cite{Hoefler_2021} for an overview of sparse deep learning.
 

Although much progress has been made in theoretical research of deep learning, it remains a challenging issue to construct an effective neural network approximation for general function spaces using as few neuron connections or neurons as possible. Most of existing network structures are specific for a particular class of functions. In this paper, we aim to propose a multi-scale sparse regularized neural network to approximate the function effectively.  A neural network with multiple hidden layers  can be viewed as a  multi-scale transformation from simple features to complex features. The layer-by-layer composite of functions can be seen as a generalization of wavelet transforms \cite{Chen2015, Daubechies1992, Micchelli_Xu1994}.
For neurons in different layers, corresponding to different transformation scales, the corresponding features have different levels of importance. Imposing different regularization parameters for different scales was proved to be an effective way to deal with multi-scale regularization problems \cite{BrennerJiangXu2009, ChenLuXuYang2008, LuShenXu2007}. Inspired by multi-scale analysis, we propose a sparse regularization network model by applying different sparse regularization penalties to the neuron connections in different layers. During the training process, the neural network adaptively learns matrix weights from given data. By sparse optimization, many weight connections are automatically zero. The remaining neural networks composed of non-zero weights form the sparse deep neural network that we desire.

This paper is organized in five sections.
In Section 2, we describe a multi-parameter regularization model for solving partial differential equations by using deep neural networks. 
We study in Section 3 the capacity of the proposed multi-parameter regularization in adaptive representing functions having certain singularities. 
In Section 4, we investigate numerical solutions of nonlinear partial differential equations by using the proposed SDNN model. Specifically, we consider two equations: the Burgers
equation and the Schr\"odinger equation since solutions of these two equations exhibit certain types of singularities. Finally, a conclusion is drawn in Section 5.

\section{A Sparse DNN Model for Solving Partial Differential Equations}

In this section, we propose a sparse DNN model for solving nonlinear partial differential equations (PDEs).

We begin with describing the PDE and its boundary, initial conditions to be considered in this paper. Suppose that $\Omega$ is an open domain in $\mathbb{R}^d$. By $\Gamma$ we denote the boundary of the domain $\Omega$. 
Let $\mathcal{F}$ denote a nonlinear differential operator, $\mathcal{I}$ the initial condition operator, and $\mathcal{B}$ the boundary operator. We consider the following boundary/initial value problem of the nonlinear partial differential equation:
\begin{align}\label{eq:PDE}
	\mathcal{F} (u(t, x)) &= 0, \ \ x \in \Omega, \ t \in [0, T],  \\
	 \label{eq:PDE:initial}
	 \mathcal{I}(u(0,x)) &=0, \ \ x \in \Omega, \ t=0, \\
	\label{eq:PDE:boundary}
	 \mathcal{B}(u(t, x))  &= 0, \ \ x \in \Gamma, \ t \in [0, T], 
\end{align}
where $T>0$, the data $u$ on $\Gamma$ and $t = 0$ are given and $u$ in $\Omega$ is the solution to be learned.
The formulation \eqref{eq:PDE}-\eqref{eq:PDE:boundary} covers a broad range of problems including conservation laws, reaction–diffusion equations, and Navier–Stokes equations. For example, the one dimensional Burgers equation can be recognized as 
$$
\mathcal{F} (u) := \ut + u \ux -  \uxx. 
$$ 


The goal of this paper is to develop a sparse DNN model for solving problem \eqref{eq:PDE}. We will conduct numerical study of the proposed model by applying it to two equations, the Burgers equation and the Schr\"odinger equation, of practical importance.

Now, we present the sparse DNN model with multi-parameter regularization.
%
%
%
%
%
We first recall the the feed forward neural network (FNN). 
A neural network can be viewed as a composition of functions.  A FNN of depth $D$ is defined to be a neural network with an input layer,  $D-1$ hidden layers, and an output layer. A neural network with more than two hidden layers is usually called a deep neural network (DNN).
Suppose that there are $d_i$ neurons in the $i$-th hidden layer. Let $W_i \in \mathbb{R}^{d_i \times d_{i-1}}$ and $b_i \in \mathbb{R}^{d_i}$ denote, respectively, the weight matrix and bias vector of the $i$-th layer. By $x_0: =x \in \mathbb{R}^{d_0}$ we denote the input vector and by $x_{i-1}\in \mathbb{R}^{d_{i-1}}$ we denote the output vector of the $(i-1)$-th layer. For the $i$-th hidden layer, we define the affine transform $L_i:\mathbb{R}^{d_{i-1}}\to \mathbb{R}^{d_i}$ by
\begin{equation*}
	L_i(x_{i-1}) := W_i x_{i-1} +b_i, \quad i = 1, 2, \dots, D.
\end{equation*} 
For an activation function $\sigma_i$, the output vector of the $i$-th hidden layer is defined as
$$ 
x_i := \sigma_i(L_i(x_{i-1})). 
$$
Given nonlinear activation functions $\sigma_i$, $i= 1, 2, \dots, D-1$, the feed forward neural network $\mathcal{N}_{\Theta}(x)$ of depth $D$ is defined as 
\begin{equation}
	\label{eq:neural_network}
\mathcal{N}_{\Theta}(x)  := L_D \circ \sigma_{D-1} \circ L_{D-1} \circ \dots \circ \sigma_{1} \circ L_1 (x),	
\end{equation}
where $\circ$ denotes the composition operator and $\Theta := \{W_i, b_i\}_{i=1}^{D}$ is the set of trainable parameters in the network. 


We first describe the physics-informed neural network (PINN) model introduced in \cite{PINN_2019} for solving the partial differential equation \eqref{eq:PDE}.
We denote by $Loss_{PDE}$ the loss of training data 
on the partial differential equation \eqref{eq:PDE}.
We choose $N_f$ collocation points $(t_f^i, x_f^i)$ by randomly sampling in domain $\Omega$ using a sampling method such as Latin hypercube sampling \cite{Helton_2003}.
We then evaluate $\mathcal{F}(\mathcal{N}_{\Theta}(t_f^i, x_f^i))$ for $i=1,2,\dots, N_f$ and define
\[
Loss_{PDE} :=\frac{1}{N_f}\sum_{i=1}^{N_f}|\mathcal{F}(\mathcal{N}_{\Theta}(t_f^i, x_f^i))|^2,
\]
where $\mathcal{F}$ is the operator for the partial differential equation \eqref{eq:PDE}.

We next describe the loss function for the boundary/initial condition.
We randomly sample $N_0$ points $x_0^i$ for the initial condition \eqref{eq:PDE:initial}, $N_b$ points $\{t_b^i, x_b^i\}$  for the boundary condition \eqref{eq:PDE:boundary}. 
The loss function  $Loss_{0}$ related to the initial value condition is given by
\[
Loss_{0} := \frac{1}{N_0}\sum_{i=1}^{N_0} |\mathcal{I}(\mathcal{N}_{\Theta}(0, x_0^i))|.
\]
The loss function $Loss_{b}$ pertaining to the boundary value is given as 
$$Loss_{b} := \frac{1}{N_{b}} \sum_{i=1}^{N_{b}} \left|\mathcal{B}(\mathcal{N}_{\Theta}(t_b^i, x_b^i) )\right|, \quad x_b^i \in \Gamma.$$

Adding the three loss functions $Loss_{PDE}$, $Loss_{0}$, and $Loss_{b}$ together gives rise to the PINN model
\begin{equation}\label{PINN}
    \min_\Theta\left\{ \frac{1}{N_f}\sum_{i=1}^{N_f}|\mathcal{F}(\mathcal{N}_{\Theta}(t_f^i, x_f^i))|^2+\frac{1}{N_0}\sum_{i=1}^{N_0} |\mathcal{I}(\mathcal{N}_{\Theta}(0, x_0^i))| +\frac{1}{N_{b}} \sum_{i=1}^{N_{b}} \left|\mathcal{B}(\mathcal{N}_{\Theta}(t_b^i, x_b^i) )\right|\right\},
\end{equation}
where $\Theta := \{W_i, b_i\}_{i=1}^{D}$.

The neural network learned from ing \eqref{PINN} is often dense and may be over-parameterized. Moreover, training data are often contaminated with noise. When noise presents, over-parameterized models may overfit training data samples and result in bad generalization to the unseen samples. 
The problem of over-fitting is often overcome by adding a regularization term:
$$
	Loss := Loss_{PDE} + \beta (Loss_{0}+ Loss_{b})+ {\rm Regularization}.
$$
The $\ell_1$- and $\ell_2$-norms are popular choices for regularization.  
Design of the regularization often makes use of prior information of the solution to be learned. It is known \cite{Candes-Romberg-Tao:CPAM:06, Donoho:IEEEIT:06, Zhang-Xu-Zhang} that the $\ell_1$-norm can promote sparsity. 
Hence, the $\ell_1$-norm regularization not only has many advantages over the $\ell_2$-norm regularization, but also leads to sparse models which can be more easily 
interpreted. Therefore, we choose to use the $\ell_1$-norm as the regularizer in this study. Furthermore, we observe that DNNs have an intrinsic multiscale structure 
whose different layers represent different scales of information, which will be validated later by numerical studies.
In fact, we will demonstrate in the next section that a smooth function or smooth parts of a function can be represented by a DNN with sparse weight matrices. This is 
because a smooth part of a function contains redundant information, which can be described very well by a few parameters, and only non-smooth parts of a function require 
more parameters to describe them. In other words, by properly choosing regularization, DNNs can lead to adaptive sparse representations of functions having certain singularities. With this understanding,
we construct an adaptive representation of a function, especially for a function having certain singularity by adopting a sparse regularization model. Our idea for the adaptive representation is to impose different sparsity penalties for different layers. Specifically, we propose a multiscale-like sparse regularization using the $\ell_1$-norm of the weight matrix for each layer with a different parameter for a different layer. The regularization with multiple parameters allows us to represent a function in a multiscale-like neural network which is determined by sparse weight matrices having different sparsity at different layers.
Such a regularization added to the loss function will enable us to robustly extract critical information of the solution of the PDE.

We now describe the proposed regularization.
For layer $i$, we denote by $W_i^{k,j}$ the $(k,j)$-th entry of matrix $W_i$, the entry in the $k$-th row and the $j$-th column. For this reason, we adopt the $\ell_1$-norm of matrix  $W_i$  defined by
$$
\|W_i\|_1 := \sum_{k=1}^{d_i} \sum_{j=1}^{d_{i-1}}|W_i^{k,j}|
$$
as our sparse regularization. Considering that different layers of the neural network play different roles in approximation of a function, we introduce here a multi-parameter regularization model
\begin{equation}
	\label{sparse:reg}
{\rm Regularization}:=	\sum_{i=1}^{D} \alpha_i \|W_i\|_1,
\end{equation}
where $\alpha_i$ are nonnegative regularization parameters. The use of different parameters for weight matrices of different layers in the regularization term  \eqref{sparse:reg} allows us to penalize the weight matrices at different layers of the neural network differently in order to extract the multiscale  representation of the solution to be learned. That is,
for a fixed $i$, parameter $\alpha_i$ determines the sparsity of weight matrix $W_i$. The larger the parameter $\alpha_i$, the more sparse the weight matrix $W_i$ is. 
The regularized loss function takes the form
\begin{equation}\label{eq:reg_inference}
		Loss :=Loss_{PDE} + \beta (Loss_{0}+ Loss_{b}) + \sum_{i=1}^{D} \alpha_i \|W_i\|_1.
\end{equation}
%
%
The parameters $\Theta := \{W_i, b_i\}_{i=1}^{D}$ of the neural network  $\mathcal{N}_{\Theta}(t, x)$ are learned by minimizing the loss function 
\begin{equation}\label{eq:minization}
	\min_{\Theta} \left\{\frac{1}{N_f}\sum_{i=1}^{N_f}|\mathcal{F}(\mathcal{N}_{\Theta}(t_f^i, x_f^i))|^2+ \beta (Loss_{0}+ Loss_{b})  + \sum_{i=1}^{D} \alpha_i \|W_i\|_1\right\}.
\end{equation}
Truncating the weights of the layers close to the input layer has an impact on all subsequent layers. In practice, we usually set smaller regularization parameters in layers close to the input and larger regularization parameters in layers close to the output. The resulting neural network will exhibit denser weight matrices near the input layer and sparser weight matrices near the output layer. This network structure reflects the multi-scale nature of neural networks and is automatically learned by sparse regularization.

Appropriate choices of the regularization parameters are key to achieve good prediction results. We need to balance sparsity and prediction accuracy. Since there are multiple regularization parameters, the regularization parameters are chosen by grid search layer by layer in this paper. In practice, we first choose the regularization parameters close to the output layer, and then gradually choose the regularization coefficients close to the input layer.

We refer equation \eqref{eq:minization} as to the sparse DNN (SDNN) model for the partial differential equation. Upon solving the minimization problem  \eqref{eq:minization}, we obtain
an approximate solution $u(t, x):= \mathcal{N}_{\Theta}(t, x)$ with sparse weight matrices. 
When the regularization parameters $\alpha_i$ are all set to 0, the SDNN model \eqref{eq:minization} reduces to the PINN model introduced in \cite{PINN_2019}. We will compare numerical performance of the proposed SDNN model with that of PINN model, for both the Burgers equation and the Schr\"odinger equation.


\section{Function Adaptive Approximation by the SDNN Model}

We explore in this section the capacity of the proposed multi-parameter regularization in adaptive representing functions that have certain singularities. We will first reveal that a DNN indeed has an intrinsic multiscale-like structure which is desirable for representing non-smooth functions. We demonstrate in our numerical studies that the proposed SDNN model can reconstruct neural networks which approximate functions in the same accuracy order with nearly the same order of network complexity, regardless the smoothness of the functions. We include in this section a numerical study of reconstruction of black holes by the proposed SDNN model. In this section, we use the rectified linear unit (ReLU) function $$
\text{ReLU}(x) := \max\{0, x\}, \ \ x\in \mathbb{R}
$$ 
as an activation function.

We first describe the data fitting problem.
Given training points $(x_i, y_i)$, $i=1,2,\dots, N$,
a non-regularized neural network is determined by minimizing the regression error, that is,
\begin{equation}
	\label{eq:non:reg}
	\min_{\Theta}  \frac{1}{N}\sum_{i=1}^{N} |\mathcal{N}_{\Theta}(x_i) - y_i|^2.
\end{equation}
The multi-parameter sparse regularization DNN model for the data fitting problem reads
\begin{equation}
	\label{eq:sparse:reg}
	\min_{\Theta} \left\{ \frac{1}{N}\sum_{i=1}^{N} |\mathcal{N}_{\Theta}(x_i) - y_i|^2 +\sum_{i=1}^{D} \alpha_i \|W_i\|_1\right\},
\end{equation}
where $\alpha_i$ are nonnegative regularization parameters and $W_i$ are weight matrices. 

In examples to be presented in this section and the section that follows, the network structure is described by the number of neurons in each layer. Specifically,  we use the notation $[d_0, d_1, \dots, d_D ]$ to describe networks that have one input layer, $D-1$ hidden layers and one output layer, with $d_0, d_1, \dots, d_D$ number of neurons, respectively. The regularization parameters, which will be presented as a vector $\alpha := [\alpha_1, \alpha_2, \dots, \alpha_D]$, are chosen so that best results are obtained. 
We will use the relative $L_2$ error to measure approximation accuracy. Suppose that $y_i$ is the exact value of function $f$ to be approximated at $x_i$, that is, $y_i = f(x_i)$, and suppose that $\hat{y}_i := \mathcal{N}_{\Theta}(x_i)$ is the output of the neural network approximation of $f$. We let $y:=[y_1, y_2,\dots, y_N]$ and $\hat{y} := [\hat{y}_1, \hat{y}_2,\dots, \hat{y}_N]$, and define the error by $\|y- \hat{y}\|_2/\|y\|_2$.  Sparsity of the weight matrices is measured by the percentage of zero entries in the weight matrices $W_i$. In our computation, we set a weight matrix entry 
$$
W_i^{k,j}=0,\ \  \mbox{if}\ \ |W_i^{k,j}|< \epsilon,
$$
where $\epsilon$ is small positive number. In our numerical examples, we set $\epsilon := 0.001$ by default.
For all numerical examples, the non-smooth, non-convex optimization problem \eqref{eq:sparse:reg} is solved by 
the Adam algorithm, which is an improved version of the stochastic gradient descent algorithm proposed in \cite{Kingma2015} for training deep learning models.



\subsection{Intrinsic adaptivity of the SDNN model}

We first investigate whether the SDNN model \eqref{eq:sparse:reg} can generate a network that has an intrinsic adaptive representation of a function. That is, a function generated by the model has a multiscale-like structure so that
the reconstructed neural networks approximate functions in the same accuracy order with nearly the same order of network complexity, regardless the smoothness of the functions. In particular, when a function is singular a sparse network with higher layers is generated to capture the higher resolution information of the function. The complexity of the network is nearly proportional to the reciprocal of the approximation error regardless whether the function is smooth or not.  
In this experiment, we consider two examples: (1) one-dimensional functions and (2) two-dimensional functions. 

In our first example, we consider approximation of the quadratic function 
\begin{equation}
	\label{fun:quadratic}
	f(x) := x^2,
\end{equation}
and the piecewise quadratic function 
\begin{equation}
	\label{fun:quadratic:discontinous}
	f(x)=\left\{
	\begin{aligned}
		&x^2+1,  \quad x \ge 0, \\
		&x^2, \quad x<0,
	\end{aligned}
	\right.
\end{equation}
by SDNN. 
Note that the function defined by \eqref{fun:quadratic} is smooth and the function by \eqref{fun:quadratic:discontinous} has a jump discontinuity at the point $0$. We applied the sparse regularized network having the architecture $[1, 10, 10, 10, 10, 1]$ to learn these functions. We divide the interval $[-2, 2]$ by the nodes $x_j:=-2+jh$, for $j:=0, 1,  \dots, 200$, with $h:=1/50$, and sample the functions $f$ at $x_j$. The test set is  $\{(x_k, f(x_k))\}$, where $x_k:=-2+kh$, $h: = 1/30$, $k= 0, 1, \dots, 120$. The network is trained by the Adam algorithm with epochs 20,000 and initial learning rate 0.001.

\begin{figure}[htbp]
	\centering
	\begin{minipage}[t]{0.48\textwidth}
		\centering
		\includegraphics[width=0.9\linewidth]{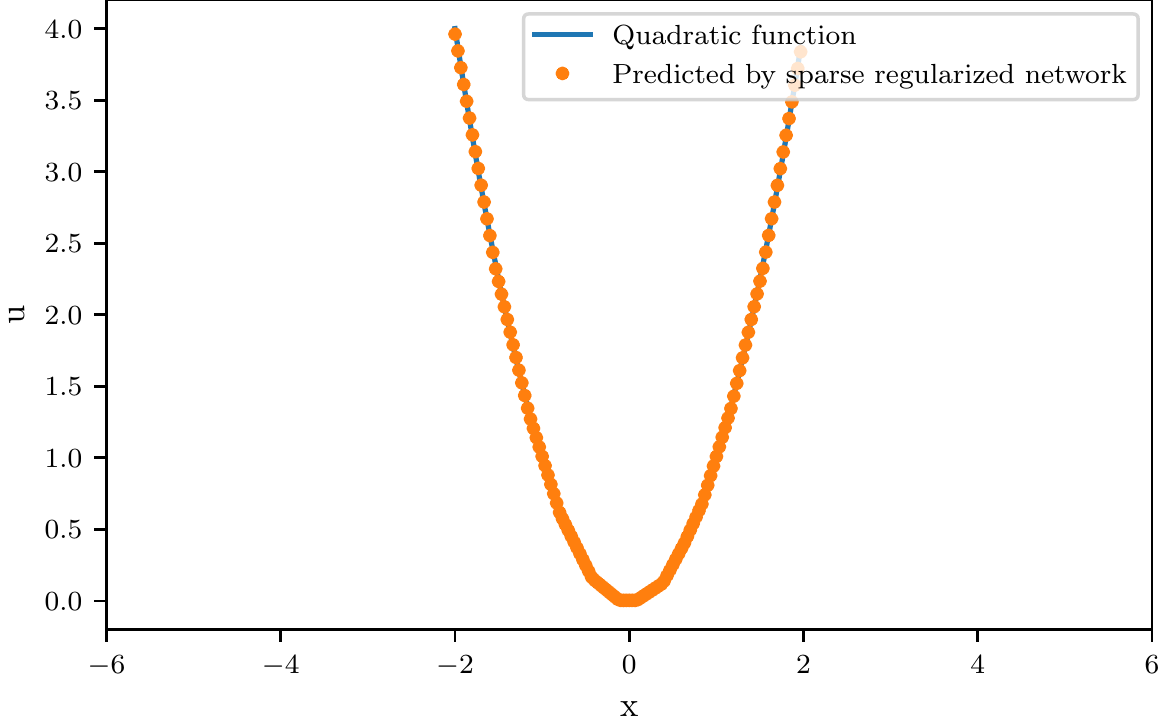}
	\end{minipage}
	\begin{minipage}[t]{0.48\textwidth}
		\centering
		\includegraphics[width=0.9\linewidth]{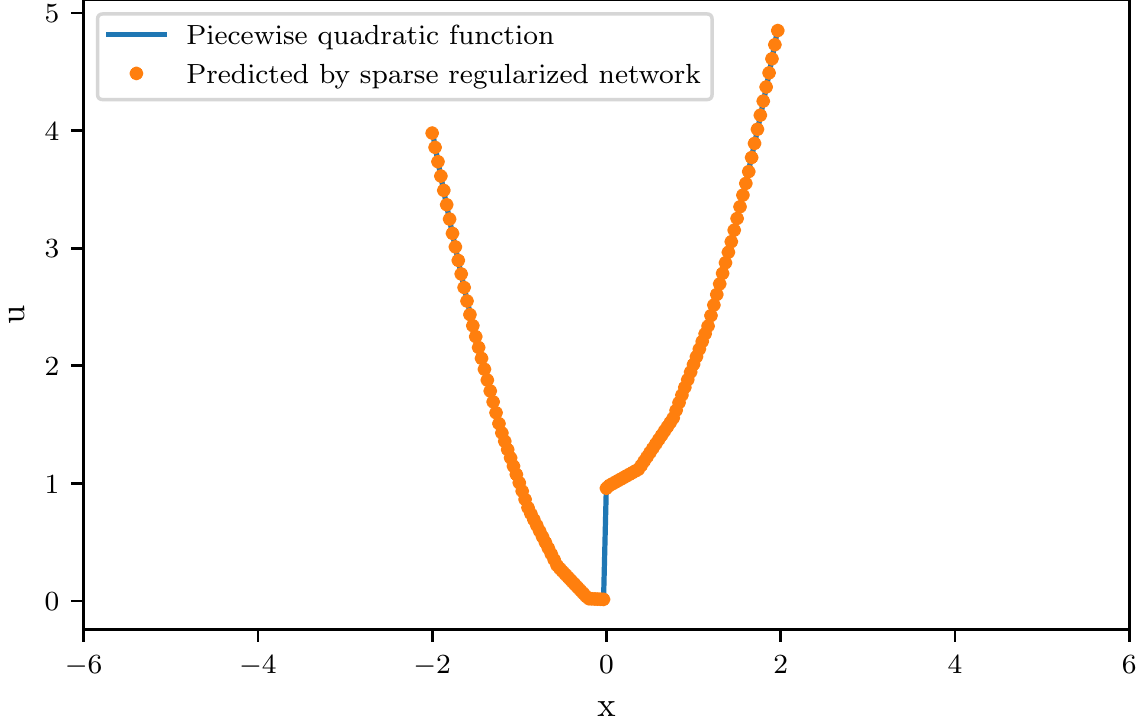}
	\end{minipage}
    \caption{Numerical results of SDNN: for function \eqref{fun:quadratic} (Left);  for function \eqref{fun:quadratic:discontinous} (Right).}
    \label{fig:quadratic}
\end{figure}

For function \eqref{fun:quadratic}, regularization parameters are set to be 
[0, 1e-4, 1e-4, 1e-3, 1e-3]. 
We obtain the prediction error 5.94e-3 for the test set. Sparsity of the resulting weight matrices is [0.0\%, 87.0\%, 95.0\%, 98.0\%, 90.0\%] and the number of nonzero weight matrix entries is 31.
Left of Figure \ref{fig:quadratic} shows the reconstructed SDNN for the function defined by  \eqref{fun:quadratic}.

\begin{table}[h]
\begin{tabular}{l|l|l}
\hline
\multirow{4}{*}{\begin{tabular}[c]{@{}l@{}}Results for \\ function \eqref{fun:quadratic}  \end{tabular}} & 
                Regularization parameters       &  [0, 1e-4, 1e-4, 1e-3, 1e-3]                 \\ \cline{2-3} 
               & Relative $L_2$ error &     5.94e-3  \\ \cline{2-3} 
                & Sparsity of weight matrices     &  [0.0\%, 87.0\%, 95.0\%, 98.0\%, 90.0\%]    \\ \cline{2-3} 
               & No. of nonzero entries      &  31                 \\ \hline
\multirow{4}{*}{\begin{tabular}[c]{@{}l@{}}Results for \\ function \eqref{fun:quadratic:discontinous}  \end{tabular}} & 
                Regularization parameters       &  [1e-5, 1e-4, 1e-4, 1e-4, 1e-3]                 \\ \cline{2-3} 
               &Relative $L_2$ error &     5.42e-3  \\ \cline{2-3} 
                & Sparsity of weight matrices     &  [50\%, 93.0\%, 88.0\%, 93.0\%, 90.0\%]     \\ \cline{2-3} 
               & No. of nonzero entries      &  32                 \\ \hline
\end{tabular}
\caption{Numerical result for quadratic function \eqref{fun:quadratic} and  piecewise quadratic function \eqref{fun:quadratic:discontinous} with network structure [1, 10, 10, 10, 10, 1].}
\label{table:quadratic_function}
\end{table}





For function \eqref{fun:quadratic:discontinous} the regularization parameters  are chosen as [1e-5, 1e-4, 1e-4, 1e-4, 1e-3]. We obtain the prediction error 5.42e-3 for the test set. Sparsity of the
resulting weight matrices is [50\%, 93.0\%, 88.0\%, 93.0\%, 90.0\%] and the number of nonzero weight matrix entries is 32. 
The reconstructed function is shown in Right of  Figure \ref{fig:quadratic}.

Numerical results for both functions  \eqref{fun:quadratic} and  \eqref{fun:quadratic:discontinous} are summarized in Table \ref{table:quadratic_function}. These results demonstrate that even though the function \eqref{fun:quadratic:discontinous} has a jump discontinuity at the point $0$, the proposed SDNN model can generate a network with nearly the same number of nonzero weight matrix entries and with the same accuracy as those for the smooth function \eqref{fun:quadratic}. This shows that the proposed SDNN model has a good adaptive approximation property.

\begin{figure}[htbp]
	\centering
	\begin{minipage}[ht]{0.5\textwidth}
		\centering
		\includegraphics[width=0.9\linewidth]{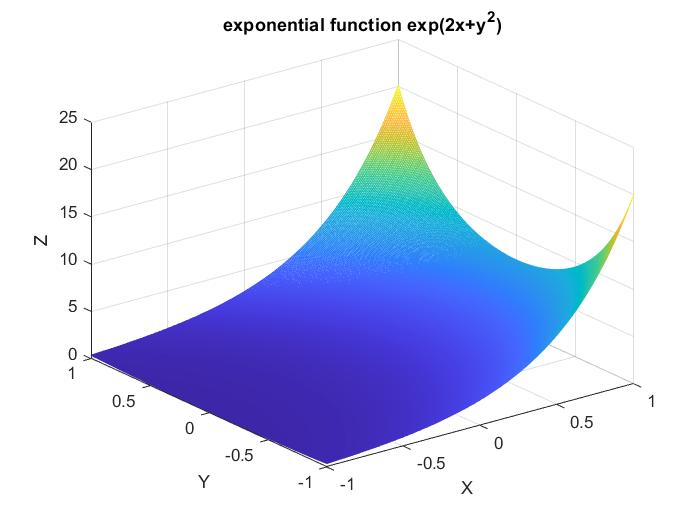}
	\end{minipage}
	\begin{minipage}[ht]{0.48\textwidth}
		\centering
		\includegraphics[width=0.9\linewidth]{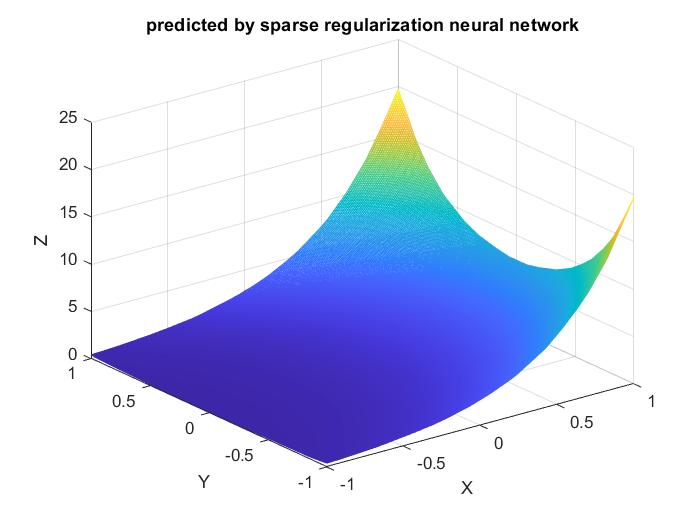}
	\end{minipage}
    \caption{Left: image of function $e^{2x+y^2}$. Right: predicted by fully connected neural network.}
    \label{fig:continuous_2d}
\end{figure}
\begin{figure}[htbp]
	\centering
	\begin{minipage}[ht]{0.48\textwidth}
		\centering
		\includegraphics[width=0.9\linewidth]{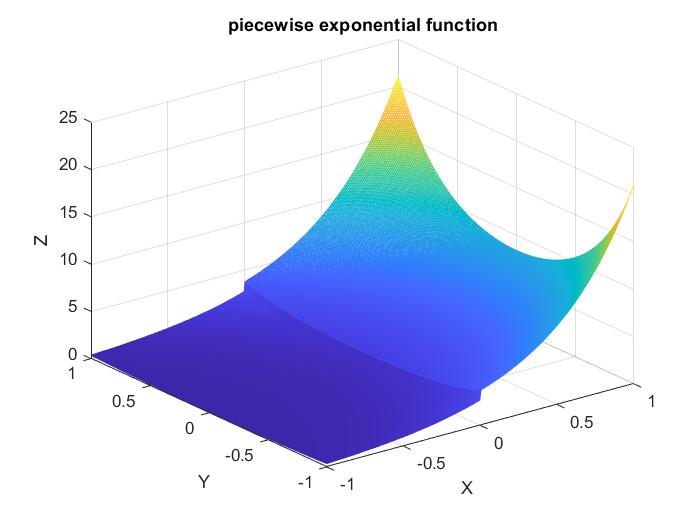}
	\end{minipage}
	\begin{minipage}[ht]{0.48\textwidth}
		\centering
		\includegraphics[width=0.9\linewidth]{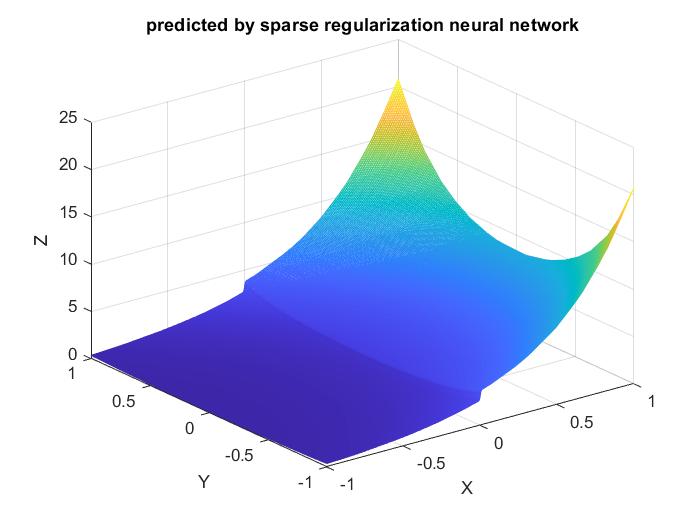}
	\end{minipage}
    \caption{Left: image of piecewise discontinuous function \eqref{eq:2d_discontinuous}. Right: predicted by sparse regularized neural network.}
    \label{fig: discontinuous_2d_predict}
\end{figure}

In our second example, we consider approximation of two-dimensional functions, once again one smooth function and one discontinuous function. We study smooth function 
\begin{equation}
	\label{eq:2d_continuous}
	g(x,y):= e^{2x+y^2},  
\end{equation}
whose image is illustrated in Figure \ref{fig:continuous_2d} (Left), and piecewise function
\begin{equation}
	\label{eq:2d_discontinuous}
	g_d(x, y)=\left\{
\begin{aligned}
	&e^{2x+y^2}+1,  \quad x \ge 0, \\
	&e^{2x+y^2}, \quad x<0,
\end{aligned}
\right.	
\end{equation}
whose image is illustrated in Figure \ref{fig: discontinuous_2d_predict} (Left). Note that function \eqref{eq:2d_continuous} is smooth and function \eqref{eq:2d_discontinuous} has a jump discontinuity along $x=0$.

For these two functions, the training data set is composed of grid points $[-1, 1]\times [-1, 1]$ uniformly discretized with step size 1/200 on $x$ and $y$ direction, and the test set is composed of grid points $[-1, 1] \times [-1, 1]$ uniformly discretized with step size 1/300 on the $x$ and $y$ directions.  The network has 2 inputs, 4 hidden layers, and 1 output, with the architecture $[2, 20, 20, 20, 20, 1]$. For each hidden layer, there are 20 neurons. The initial learning rate for Adam is set to 0.001. The batch size is equal to 1024.%

For function \eqref{eq:2d_continuous} we set the sparse regularization parameters as [0, 1e-6, 1e-4, 1e-4, 1e-4]. After 10,000 epochs training, the sparsity of weight matrices is [0.0\%, 68.5\%, 95.75\%, 97.75\%, 80.0\%] and the number
of nonzero weight matrix entries is 178. The prediction error for the test set is 4.38e-3. For function \eqref{eq:2d_discontinuous}, the regularization parameters are set to be [1e-4, 1e-5, 1e-5, 1e-4, 1e-4]. 
The sparsity of weight matrices after regularization are [60.0\%, 71.75\%, 81.75\%, 97.5\%, 90.0\%] and the number
of nonzero weight matrix entries is 206. 
The prediction error for sparse regularized deep neural network is  4.27e-3, which is even slightly better than that for function \eqref{eq:2d_continuous}.
The images of the reconstructed functions are shown respectively in Figures \ref{fig:continuous_2d}, \ref{fig: discontinuous_2d_predict} (Right). Numerical results for this example are reported in Table \ref{table:2d:function}.


\begin{table}[h]
\begin{tabular}{l|l|l}
\hline
\multirow{4}{*}{\begin{tabular}[c]{@{}l@{}}Results for \\ function \eqref{eq:2d_continuous}  \end{tabular}} & 
                Regularization parameters       &  [0, 1e-6, 1e-4, 1e-4, 1e-4]                 \\ \cline{2-3} 
               & Relative $L_2$ error &     4.38e-3  \\ \cline{2-3} 
                & Sparsity of weight matrices     &   [0.0\%, 68.5\%, 95.75\%, 97.75\%, 80.0\%]    \\ \cline{2-3} 
               & No. of nonzero connections      &   178               \\ \hline
\multirow{4}{*}{\begin{tabular}[c]{@{}l@{}}Results for \\ function \eqref{eq:2d_discontinuous}  \end{tabular}} & 
                Regularization parameters       &  [1e-4, 1e-5, 1e-5, 1e-4, 1e-4]                 \\ \cline{2-3} 
               &Relative $L_2$ error &     4.27e-3  \\ \cline{2-3} 
                & Sparsity of weight matrices     &   [60.0\%, 71.75\%, 81.75\%, 97.5\%, 90.0\%]    \\ \cline{2-3} 
               & No. of nonzero connections      &  206                \\ \hline
\end{tabular}
\caption{Numerical result for two dimensional function \eqref{eq:2d_continuous} and \eqref{eq:2d_discontinuous} with network structure [2, 20, 20, 20, 20, 1].} 
\label{table:2d:function}
\end{table}

The numerical results presented in this subsection indicate that indeed the proposed SDNN model has an excellent adaptivity property in the sense that it generates networks with nearly the same number of nonzero weight matrix entries and the same order of approximation accuracy for functions regardless their smoothness.

\subsection{An example of adaptive function approximation by the SDNN model}

The second experiment is designed to test the sparsity of the network learned from the SDNN model \eqref{eq:sparse:reg} and the model's generalization ability. Specifically, in this example, we demonstrate that the sparse model \eqref{eq:sparse:reg} leads to a sparse DNN with higher accuracy in comparison to the standard DNN model \eqref{eq:non:reg}. We consider the absolute value function 
$$
y=f(x):=|x|, \ \ \mbox{for}\ \ x\in \mathbb{R}.
$$
Note that function $f$ is not differentiable at $x=0$.

We adopt the same network architecture, that is, 1 input layer, 2 hidden layers, 1 output layer  and each hidden layer containing 5 neurons, for both the standard DNN model \eqref{eq:non:reg} and the SDNN model \eqref{eq:sparse:reg}. The training set is composed of equal-distance  grid points laying in $[-2, 2]$ with step size 0.01. The test set is composed of equal-distance grid points in $[-5, 5]$ with step size 0.1. For the sparse regularized network, the regularization parameters are set as [1e-4, 1e-3, 1e-3].
For both the standard DNN model and the  SDNN model, the number of epoch equals 10,000. The initial learning rate is set to 0.001.



\begin{table}[h]
	\centering 
\begin{tabular}{c||c||c}
	\hline
	& \multirow{2}{*}{Relative $L_2$ error}  & Sparsity of weight matrices \\
	&               & [$W_1$, $W_2$, $W_3$] \\
	\hline
	Standard DNN model & 5.58e-2 & [0\%, 0\%, 0\%] \\
	\hline
	SDNN model &  1.87e-3 & [20\%, 92\%, 80\%] \\
	\hline
\end{tabular}	
\caption{Approximation of the absolute value function by a SDNN with regularization parameters  [1e-4, 1e-3, 1e-3].} 
\label{table:absolute_function}	
\end{table}

We present numerical results of this experiment in Table \ref{table:absolute_function}, where we compare errors and sparsity of the functions learned from the two models.
Clearly, the network learned from the standard DNN model 
is non-sparse: all entries of its weight matrices are nonzero. While the network learned from the SDNN model has a good sparsity property: There are only 1 non-zero entries in $W_3$ and 2 non-zero entries in $W_2$ in the network learned from the SDNN model. 
Note that the absolution value function is the linear composition of two ReLU functions, that is $|x| = \text{ReLU}(x) + \text{ReLU}(-x)$. The SDNN model is able to find a linear combination of the two functions to represent the function $f(x):=|x|$ but the standard DNN model fails to do so.

\begin{figure}[htbp]
	\centering

	\begin{minipage}[t]{0.45\textwidth}
		\centering
		\includegraphics[width=0.9\linewidth]{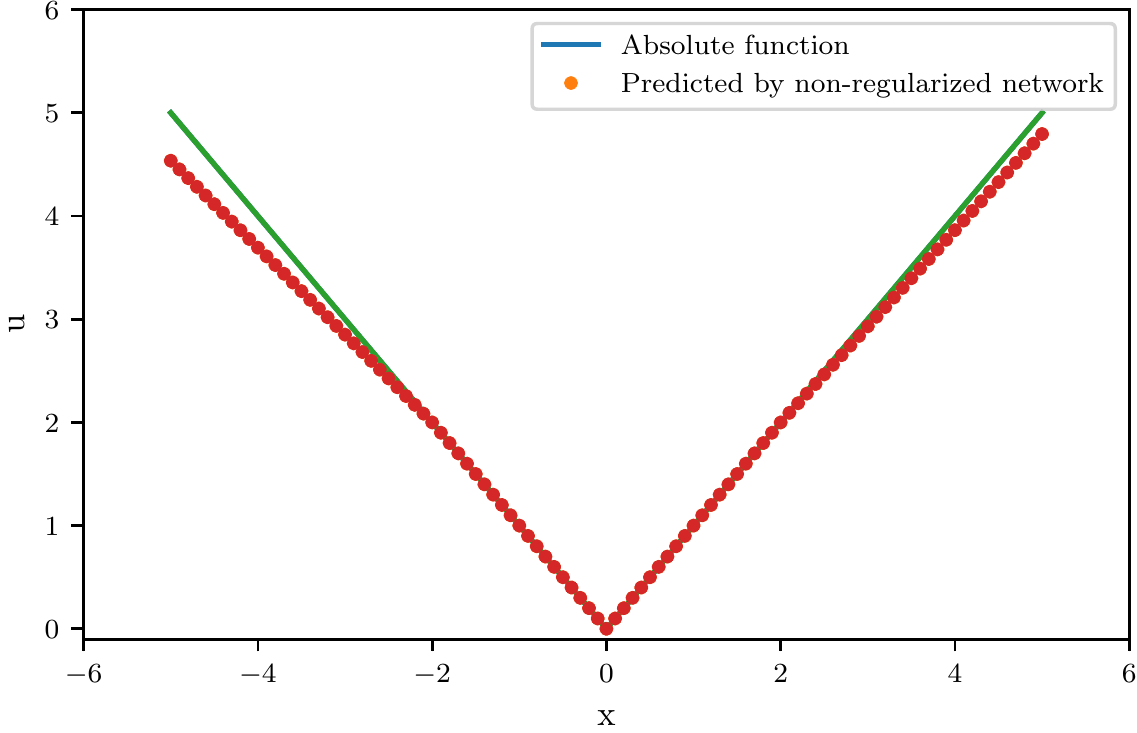}
		\caption{Reconstruction of function $f(x):=|x|$ by the standard DNN model \eqref{eq:non:reg}.}
		\label{fig:absolutefun3layern}
	\end{minipage}
	\quad \,
	\begin{minipage}[t]{0.45\textwidth}
	\centering
	\includegraphics[width=0.9\linewidth]{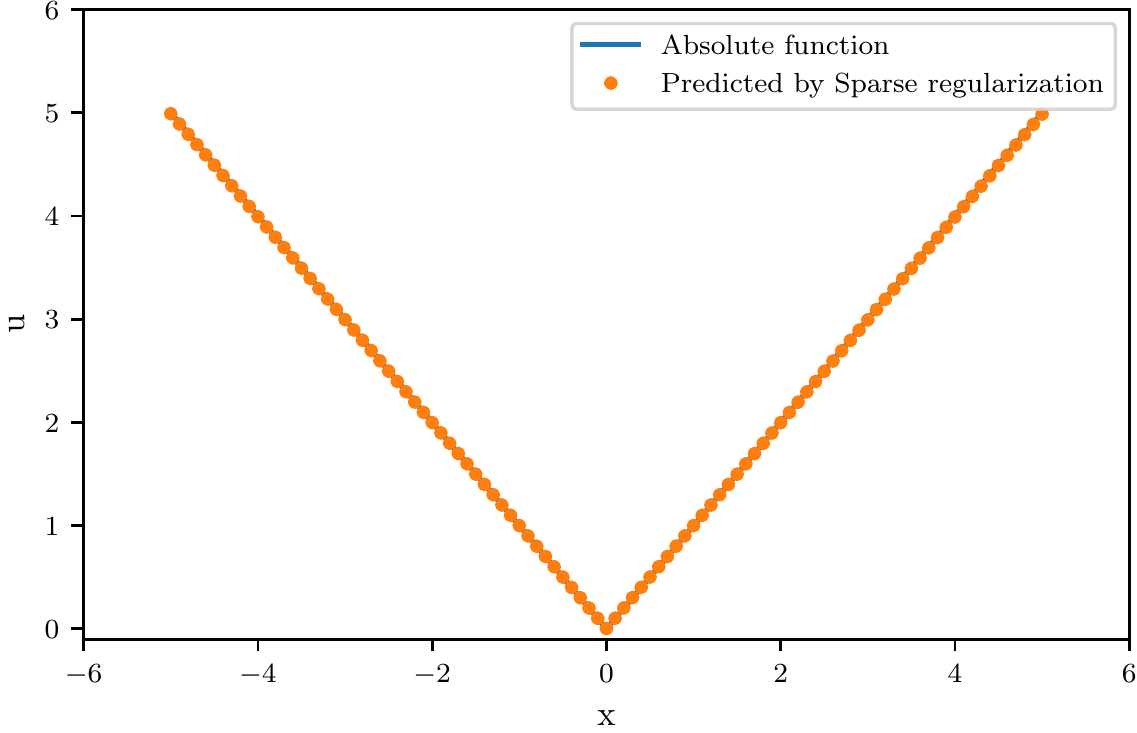}
	\caption{Reconstruction of function $f(x):=|x|$ by the SDNN model  \eqref{eq:sparse:reg}.}
	\label{fig:absolutefun3layer10000}
    \end{minipage}
\end{figure}

We plot the graphs of the reconstructed functions by the standard DNN model \eqref{eq:non:reg} and  the SDNN model  \eqref{eq:sparse:reg} in Figures \ref{fig:absolutefun3layern} and \ref{fig:absolutefun3layer10000}, respectively.
It can be seen from  Figure \ref{fig:absolutefun3layern} that the function 
reconstructed by the standard DNN model \eqref{eq:non:reg} has large errors in 
the interval $[3, 5]$. Figure \ref{fig:absolutefun3layer10000} shows that the function reconstructed by the SDNN model  \eqref{eq:sparse:reg} almost coincides with the original function. This example indicates that the SDNN model  \eqref{eq:sparse:reg}  has better generalization ability than the standard DNN model \eqref{eq:non:reg}.

\subsection{Reconstruction of a black hole}

\begin{figure}[htbp]
	\centering
	\begin{minipage}[ht]{0.4\textwidth}
		\centering
		\includegraphics[width=0.9\linewidth]{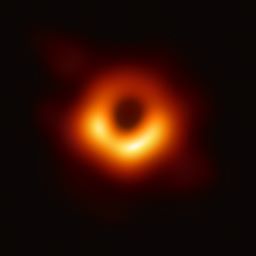}
	\end{minipage}
	\begin{minipage}[ht]{0.4\textwidth}
		\centering
		\includegraphics[width=0.9\linewidth]{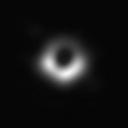}
	\end{minipage}
    \caption{Left: color image of the black hole. Right: gray image of the black hole.}
     \label{black_hole}
\end{figure}

In this example, we consider reconstruction of the image of a black hole by the SDNN model. Specifically, we compare numerical results and reconstructed image quality of the SDNN model with those of the standard DNN model. We choose a color image of the black hole shown in Figure \ref{black_hole} (Left), which is turned into a gray image shown in Figure \ref{black_hole} (Right).
The image has the size $128 \times 128$ and can be represented as a two-dimensional discrete function. The value of the gray image at the point $(x_1, x_2)$ is defined as a  $f_{image}(x_1, x_2)$, $x_1, x_2 =1, 2, \dots, 128$. The function clearly has singularities. 

\begin{figure}[!htbp]
	\centering
	\begin{minipage}[ht]{0.4\textwidth}
		\centering
		\includegraphics[width=0.9\linewidth]{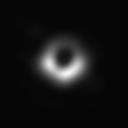}
	\end{minipage}
	\begin{minipage}[ht]{0.4\textwidth}
		\centering
		\includegraphics[width=0.9\linewidth]{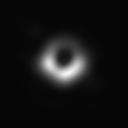}
	\end{minipage}
    \caption{Images of the black hole reconstructed: by the standard DNN model (Left) and by the SDNN model (Right).}
    \label{result:blackhole:fnn:sr}
\end{figure}

The network architecture that we used for the construction is 
$$
[2, 100, 100, 100, 100, 100, 100, 1].
$$
We randomly choose 5,000 points $(x_1^{i}, x_2^{i}, f_{image}(x_1^{i}, x_2^{i}))$ by uniform sampling, $i=1, 2, \dots, 5,000$, from the image of the black hole to train both the standard neural network and the sparse regularized network. 
The optimizer is chosen as the Adam algorithm with batch size 1,024. The number of epoch is 40,000. The patience parameter of early stopping is 200.
Prediction results by the standard DNN model and by the SDNN model are shown respectively on Figure \ref{result:blackhole:fnn:sr} (Left) and (Right).

\begin{table}[h]
	\centering 
\begin{tabular}{l||l}
	\hline
	\multicolumn{2}{c}{Standard DNN model}\\
	\hline
	 Relative $L_2$ error & Sparsity of weight matrices  \\
	\hline
	 9.66e-3              & [0\%, 0\%, 0\%,  0\%,  0\%, 0\%, 0\%] \\
	\hline
	\multicolumn{2}{c}{SDNN model} \\
	\hline
	Relative $L_2$ error & Sparsity of weight matrices  \\
	\hline 
	 9.28e-3  & [44.0\%, 78.3\%, 78.4\%, 80.3\%, 96.5\%, 98.2\%, 84.0\%]
 \\
	\hline
\end{tabular}	
\caption{Numerical results of the black hole reconstructed by the standard DNN model vs. the SDNN model.} 
\label{table:blackhole}	
\end{table}

\begin{figure}[!htbp]
	\centering
	\begin{minipage}[ht]{0.48\textwidth}
		\centering
		\includegraphics[width=0.9\linewidth]{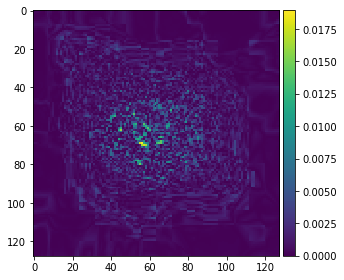}
	\end{minipage}
	\begin{minipage}[ht]{0.48\textwidth}
		\centering
		\includegraphics[width=0.9\linewidth]{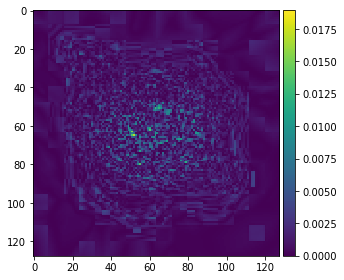}
	\end{minipage}
    \caption{Reconstruction errors of the black hole: by the standard DNN model (Left) and by the SDNN model (Right).}
    \label{error:blackhole:fnn:sr}
\end{figure}

Error images of the two models  are presented in Figure \ref{error:blackhole:fnn:sr}, from which
it can be seen that the sparse network has a smaller reconstruction error. 
The prediction error of the fully connected network is 9.66e-3. For the sparse regularized neural network, the regularized parameters are set to be [1e-9, 1e-9, 1e-9, 1e-9, 1e-8, 1e-8, 1e-8].
The prediction error of the sparse regularized network is 9.28e-3. 
The sparsity of the weight matrices are
[44.0\%, 78.3\%, 78.4\%, 80.3\%, 96.5\%, 98.2\%, 84.0\%].
It shows that the sparse regularized network uses fewer neurons and has smaller prediction error.
This indicates that by using the proposed multi-parameter sparse regularization, the deep neural network has the ability of multi-scale and adaptive learning.


\section{Numerical Solutions of Partial Differential Equations}

We study in this section numerical performance of the proposed SDNN model for solving partial differential equations. We consider two equations: the Burgers equation and the Schr\"odinger equation. For both of these two equations, we choose the hyperbolic tangent (tanh) function defined by
$$
\tanh(x):=\frac{e^x – e^{-x}}{e^x + e^{-x}}, \ \ x\in\mathbb{R}
$$
as the activation function to build networks for our approximate solutions due to its differentiability which is required by the differential equations. 

\begin{figure}[!h]
	\includegraphics[width = 1.0\textwidth]{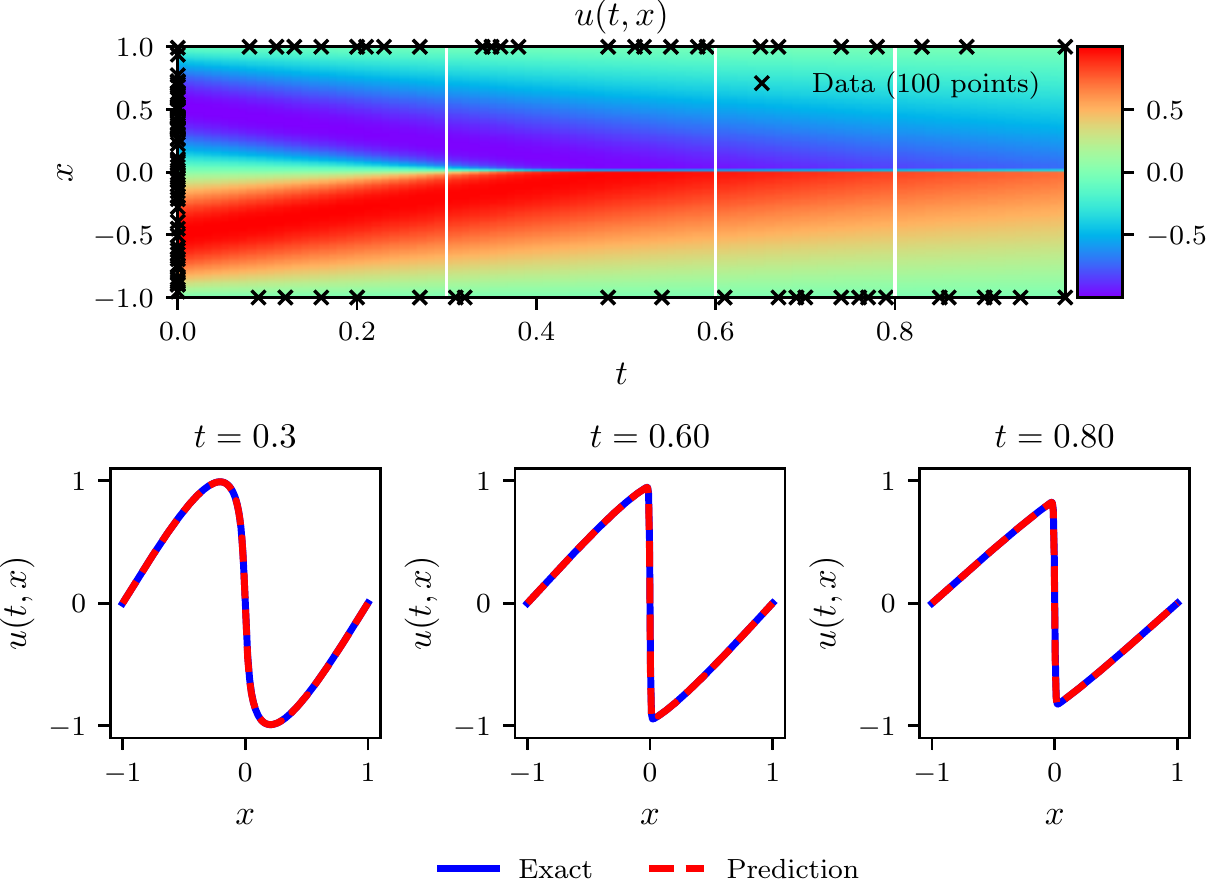}
	\vspace*{-1.0cm}
	\caption{{\em Burgers equation:} {\it Top:} The training data and predicted solution $u(t, x)$ for sparse deep neural network with [2, 50, 50, 50, 1],  regularization parameter $\alpha = \text{[1e-6, 1e-6, 1e-6, 1e-4]}$, $\beta = 20$. {\it Bottom:} Predicted solution $u(t, x)$ at time $t=0.3$, $t=0.6$, and $t=0.8$.}
	\label{fig:Burgers_CT_identification}
\end{figure}

\subsection{The Burgers equation}
The Burgers equation has attracted much attention since it is often used as simplified model for turbulence and shock waves \cite{Liu1995}. It is well-known that the solution of this equation presents a jump discontinuity (a shock wave), even though the initial function is smooth. 

In this example, we consider the following one dimensional Burgers equation
\begin{eqnarray}\label{eq:Burgers}
	&&u_t (t,x)+  u(t,x) u_x(t,x) - \frac{0.01}{\pi}  u_{xx} (t,x)= 0, \ t \in (0,1],  x\in (-1, 1), \\
	\label{eq:Burgers:initial}
	&&u(0,x) = -\sin(\pi x),  \\
	\label{eq:Burgers:boundary}
	&&u(t,-1) = u(t, 1) = 0.  
\end{eqnarray}
The analytic solution of this equation, known in \cite{Basdevant1986}, will be used as our exact solution for comparison.
Indeed, the analytic solution has the form 
$$
u(t, x):=-\frac{\int_{-\infty}^{+\infty} \sin \pi(x-\eta) h(x-\eta) \exp \left(-\eta^{2} / 4 \nu t\right) \mathrm{d} \eta }
{\int_{-\infty}^{+\infty} h(x-\eta) \exp \left(-\eta^{2} / 4 \nu t\right) \mathrm{d} \eta},\ \ t \in [0,1], \ x\in [-1, 1],
$$
where $\nu := 0.01/\pi$ and $h(y):=\exp (-\cos \pi y / 2 \pi \nu)$.
A neural network solution of equation \eqref{eq:Burgers}-\eqref{eq:Burgers:boundary} was obtained recently from the standard DNN model in \cite{PINN_2019}.

We apply the setting  \eqref{eq:PDE}-\eqref{eq:PDE:boundary}  with 
$$
\mathcal{F}(u(t, x)) := u_t(t,x) +  u(t,x) u_x(t,x) -\frac{0.01}{\pi}  u_{xx} (t,x),\  t \in (0,1], \ x\in (-1, 1).
$$
Let $\left\{x_{0}^{i}, u_0^{i}\right\}_{i=1}^{N_{0}}$ denote the training data of $u$ satisfying initial condition \eqref{eq:Burgers:initial}, that is, $u_0^{i} = -\sin(\pi x_{0}^{i})$. Let $\left\{ t_{b_1}^{i} \right\}_{i=1}^{N_{b_1}}$ and $\left\{ t_{b_2}^{i} \right\}_{i=1}^{N_{b_2}}$ be the collocation points related to boundary condition \eqref{eq:Burgers:boundary} for $x = -1$ and $x = 1$ respectively. We denote by $\left\{t_{f}^{i}, x_{f}^{i}\right\}_{i=1}^{N_{f}}$ the collocation points for $\mathcal{F}(u(t, x))$ in $[0, 1] \times (-1, 1)$.
The sparse deep neural network $\mathcal{N}_{\Theta}(t, x)$ are learned by model \eqref{eq:minization} with
%
$$
   Loss_{0} =   \frac{1}{N_0}\sum_{i=1}^{N_0} |\mathcal{N}_{\Theta}(0, x_0^i) - u_0^i|^2,
$$
and
$$
   Loss_{b} = \frac{1}{N_{b_1}} \sum_{i=1}^{N_{b_1}} \left|\mathcal{N}_{\Theta}(t_{b_1}^i, -1) \right| + \frac{1}{N_{b_2}}\sum_{i=1}^{N_{b_2}} \left|\mathcal{N}_{\Theta}(t_{b_2}^i, 1) \right|.
$$


\begin{table}[!h]
	\centering
	\begin{tabular}{c||c||c}
		\hline
	    \multirow{2}{*}{Algorithms} &Parameters $\alpha$  & Relative      \\
	     & \& sparsity of weight matrices  & $L_2$ error \\
	     \hline
		\multirow{2}{*}{PINN}  & No regularization  & \multirow{2}{*}{2.45e-2} \\
		    & [0.0\%, 0.2\%, 0.5\%, 0.0\%]  &  \\
		\hline
	    SDNN  	 
		      &[1e-6, 1e-6, 1e-6, 1e-4]	
		      & \multirow{2}{*}{1.68e-3}  \\ 
			    ($\beta =20$)       &   [13.0\%, 61.9\%, 71.2\%, 62.0\%]                        & \\ \hline
	\end{tabular}
	\caption{{\em The Burgers equation:} A neural network of 4 layers, with network architecture [2, 50, 50, 50, 1]. } \label{tab:Burgers_inference}
\end{table}

In this experiment, 100 data points are randomly selected from boundary and initial data points, among which $N_{b_1} = 25$ points are located on the boundary $x=-1$, $N_{b_2} = 23$ points on the boundary $x=1$, and $N_0 = 52$ points on the initial line $t=0$. The distribution of random collocation points is shown in the top of Figure \ref{fig:Burgers_CT_identification}.
The number of collocation points of the partial differential equation is $N_f =10,000$ by employing the Latin hypercube sampling method. 
The test set is composed of grid points $[0, 1] \times [-1, 1]$ uniformly discretized with step size 1/100 on the $t$ direction  and step size 2/255 on the $x$ direction.

\begin{table}[!h]
	\centering
	\begin{tabular}{c||c||c}
		\hline
	    \multirow{2}{*}{Algorithms} &Parameters $\alpha$  & Relative      \\
	     & \& sparsity of weight matrices  & $L_2$ error \\
	     \hline
		\multirow{2}{*}{PINN}  & No regularization  & \multirow{2}{*}{3.39e-3} \\
		    & [0.0\%, 0.8\%, 0.6\%, 0.6\%, 0.8\%, 0.6\%, 0.4\%, 0.0]  &  \\
		\hline
	   SDNN	        
		& [0, 0, 0, 0, 1e-7, 1e-10, 1e-6, 1e-5] &  \multirow{2}{*}{1.45e-4} \\
	  ($\beta =10$)&  [0.0\%, 0.7\%, 0.8\%, 0.7\%, 15.6\%, 0.5\%, 93.8\%, 94.0\%]                         & \\ \hline
	   SDNN	 & [1e-6, 1e-6, 1e-6, 1e-6, 1e-6, 1e-6, 1e-5, 1e-5] &  \multirow{2}{*}{4.83e-4} \\
	  ($\beta =10$)&  [25.0\%, 78.6\%, 85.3\%, 82.8\%, 79.5\%, 84.0\%, 98.6\%, 94.0\%]                         & \\ \hline		
	\end{tabular}
	\caption{{\em The Burgers equation:} Neural networks of 8 layers with  network architecture [2, 50, 50, 50, 50, 50, 50, 50, 1]. } \label{tab:Burgers_inference_8}
\end{table}

We use two different network architectures [2, 50, 50, 50, 1] and [2, 50, 50, 50, 50, 50, 50, 50, 1] for DNNs. We choose Adam as the optimizer for both neural networks. The number of epoch is 30,000. The initial learning rate is set to 0.001.
Numerical results of these two networks presented respectively in Tables \ref{tab:Burgers_inference} and \ref{tab:Burgers_inference_8} show that the proposed SDNN model outperforms the PINN model in both weight matrix sparsity and approximation accuracy.



\subsection{The Schr\"odinger equation}

The Schr\"odinger equation is the most essential equation of non-relativistic quantum mechanics. It plays an important role in studying nonlinear optics, Bose-Einstein condensates, protein folding and bending. It is also a model equation for studying waves propagation and soliton \cite{Serkin2000}.
In this subsection, we consider a one-dimensional Schr\"odinger equation with periodic boundary conditions
\begin{equation}
	\label{eq:Sch}
\begin{aligned}
	&i u_{t}(t,x)+0.5 u_{x x}(t,x)+|u(t,x)|^{2} u(t,x)=0, \ t \in (0, \pi / 2] , \ x \in(-5, 5), \\
	&u(0, x)=2 \operatorname{sech}(x), \\
	&u(t, -5)=u(t, 5), \\
	&u_{x}(t, -5)=u_{x}(t, 5).
\end{aligned}
\end{equation}
Note that the solution $u$ of problem \eqref{eq:Sch} is a complex-valued function. The goal of this study is to test the effectiveness of the proposed SDNN model in solving complex-valued nonlinear differential equations with periodic boundary conditions, with a comparison to the standard DNN model recently developed in \cite{PINN_2019}.

Problem \eqref{eq:Sch} falls into the setting \eqref{eq:PDE}-\eqref{eq:PDE:boundary} with 
$$
\mathcal{F}(u(t, x)) := i u_{t}(t,x)+0.5 u_{x x}(t,x)+|u(t,x)|^{2}u(t,x),\ t \in (0, \pi / 2] , \ x \in(-5, 5).
$$
Let $\psi$ and $\phi$ be respectively the real part and imaginary part of the solution $u$ of problem \eqref{eq:Sch}. 
We intend to approximate the solution $u(t,x)$ by a neural network $\mathcal{N}_{\Theta}(t, x)$ with two inputs $(t, x)$ and two outputs which approximate $\psi(t, x)$ and $\phi(t, x))$, respectively. 
Let $\left\{x_{0}^{i}, u_{0}^{i}\right\}_{i=1}^{N_{0}}$ denote the training data to enforce the initial condition at time $t=0$, that is, $u_{0}^{i} = \operatorname{sech}(x_{0}^{i})$, $\left\{t_{b}^{i}\right\}_{i=1}^{N_{b}}$ the collocation points on the boundary $x=-5$ and $x=5$ to enforce the periodic boundary conditions, and $\left\{t_{f}^{i}, x_{f}^{i}\right\}_{i=1}^{N_{f}}$ the collocation points in $ (0, \pi/2] \times (-5, 5)$. These collocation points were generated by the Latin hypercube sampling method.
We then learn the neural network  $\mathcal{N}_{\Theta}(t,x)$ by  model \eqref{eq:minization} with
$$
  Loss_{0} :=    \frac{1}{N_{0}} \sum_{i=1}^{N_{0}}\left|\mathcal{N}_{\Theta}\left(0, x_{0}^{i}\right)-u_{0}^{i}\right|^{2},
$$
and
$$
Loss_{b} : =\frac{1}{N_{b}} \sum_{i=1}^{N_{b}}\left(\left|\mathcal{N}_{\Theta}\left(t_{b}^{i}, -5\right)-\mathcal{N}_{\Theta}\left(t_{b}^{i}, 5\right)\right|^{2}+\left|\frac{\partial\mathcal{N}_{\Theta}}{\partial x}\left(t_{b}^{i}, -5\right)-\frac{\partial\mathcal{N}_{\Theta}}{\partial x}\left(t_{b}^{i}, 5\right)\right|^{2}\right).
$$


\begin{figure}[!htb]
	\includegraphics[width = 1.0\textwidth]{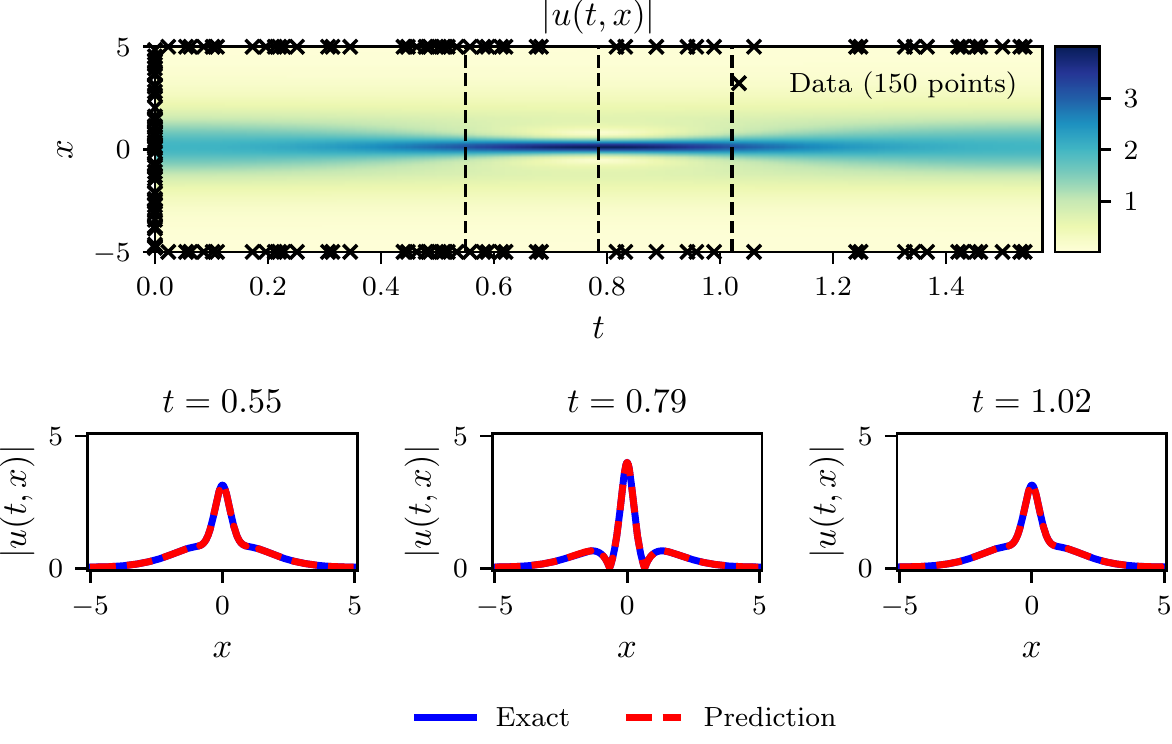}
	\vspace*{-1.0cm}
	\caption{{\em The Schrödinger equation:} {\it Top:} The training data and predicted solution $|u(t,x)|$ by SDNN with network architecture [2,
50, 50, 50, 50, 50, 50, 2],  regularization parameters $\alpha := \text{[9e-7, 5e-7, 6e-7, 7e-7, 8e-7, 1e-6, 1e-5]}$, and $\beta := 10$. {\it Bottom:} Predicted solutions at time $t:=0.55$, $t:=0.79$, and $t:=1.02$.}
	\label{fig:Sch}
\end{figure}

A reference solution of  problem \eqref{eq:Sch}  is solved by a Fourier spectral method using the Chebfun
package \cite{Driscoll:2014}. Specifically, we obtain the reference solution by using 256 Fourier modes for space discretization and an explicit fourth-order Runge–Kutta method (RK4) with time-step $\Delta t := (\pi/2) \times 10^{-6}$ for time discretization. 
For more details of the discretization of Schr\"odinger equation \eqref{eq:Sch}, the readers are referred to  \cite{PINN_2019}.

For both the standard network and the sparse network, we used the network architecture [2, 50, 50, 50, 50, 50, 50, 2]. Both the networks were trained by the Adam algorithm with 30,000 epochs. The initial learning rate is set to 0.001.
The training set is composed of $N_0 := 50$ data points on $u(0, x)$,  $N_b := 50$ sample points for enforcing the periodic boundaries, and $N_f := 20,000$ sample points inside the solution domain of equation \eqref{eq:Sch}. The test set is composed of grid points $(0, \pi/2] \times [-5, 5]$ uniformly discretized with step size $\pi$/400 on the $t$ direction and step size 10/256 on the $x$ direction.

Numerical results for this example are listed in Table \ref{tab:Sch_8}. 
As we can see, the sparse network has a smaller prediction error than the standard network. When regularization parameters $\alpha$ = [0, 0, 0, 0, 5e-7, 1e-6, 1e-5], the relative $L_2$ error is smaller than the PINN method. When regularization parameters $\alpha$ are taken as [9e-7, 5e-7, 6e-7, 7e-7, 8e-7, 1e-6, 1e-5], the sparsity of weight matrices are [22.0\%, 50.5\%, 51.9\%, 50.6\%, 50.0\%, 64.5\%, 66.0\%]. In other words, after removing more than half of the neural network connections, the sparse neural network still has a slightly higher prediction accuracy. The predicted solution of the SDNN is illustrated in Figure \ref{fig:Sch}. These numerical results clearly confirm that the proposed SDNN model outperforms the standard DNN model.

\begin{table}[!h]
	\centering
	\begin{tabular}{c||c||c}
		\hline
	    \multirow{2}{*}{Algorithms} &Parameters $\alpha$  & Relative      \\
	     & \& sparsity of weight matrices  & $L_2$ error \\
	     \hline
		\multirow{2}{*}{PINN}  & No regularization  & \multirow{2}{*}{1.41e-3} \\
		      & [0.0\%, 0.3\%, 0.4\%, 0.6\%, 0.4\%, 0.68\%, 0.0\%]                    &     \\     
		\hline
		SDNN  & $\alpha$ =[0, 0, 0, 0, 5e-7, 1e-6, 1e-5]  & \multirow{2}{*}{8.15e-4} \\
		 ($\beta =10) $ & [0.7\%, 0.3\%, 0.4\%, 50.6\%, 38.0\%, 74.7\%, 77.0\%]  &   \\
		 	\hline
		 SDNN  & $\alpha$ =[9e-7, 5e-7, 6e-7, 7e-7, 8e-7, 1e-6, 1e-5]  & \multirow{2}{*}{1.38e-3} \\
		 ($\beta =10) $ & [22.0\%, 50.5\%, 51.9\%, 50.6\%, 50.0\%, 64.5\%, 66.0\%]  &   \\
		\hline
	\end{tabular}
	\caption{{\em The Schrödinger equation:}  The neural network of 7 layers with network architecture [2, 50, 50, 50, 50, 50, 50, 2]. } 
	\label{tab:Sch_8}
\end{table}

\section{Conclusion}
A sparse network requires less memory and computing time to operate it and thus it is desirable.
We have developed a sparse deep neural network model by employing a sparse regularization with  multiple parameters for solving nonlinear partial differential equations. Noticing that neural networks are layer-by-layer composite structures with an intrinsic multi-scale structure, we observe that the network weights of different layers have different weights of importance. Aiming at generating a sparse network structure while maintaining approximation accuracy, we proposed to impose different regularization parameters on different layers of the neural network. We first tested the proposed sparse regularization model in approximation of singular functions, and discovered that the proposed model can not only generate an adaptive approximation of functions having singularities but also have better generalization than the standard network. We then developed a sparse deep neural network model for solving nonlinear partial differential equations whose solutions may have certain singularities. Numerical examples show that the proposed model can remove redundant network connections leading to sparse networks and has better generalization ability. Theoretical investigation will be performed in a follow-up paper.


\begin{thebibliography}{100}
    \bibitem{Adcock2021}
	B. Adcock and N. Dexter, The gap between theory and practice in function approximation with deep neural networks, SIAM Journal on Mathematics of Data Science, 3 (2) (2021): 624-655.
	
	\bibitem{Basdevant1986}
	C. Basdevant, M. Deville, P. Haldenwang, J. Lacroix, J. Ouazzani, R. Peyret, P. Orlandi, and A. Patera, Spectral and finite difference solutions of the Burgers equation, Comput. Fluids, 14 (1986): 23-41.

   \bibitem{Bolcskei2019}
   H. Bölcskei, P. Grohs, G. Kutyniok, and P. Petersen, Optimal approximation with sparsely connected deep neural networks, SIAM Journal on Mathematics of Data Science, 1 (1) (2019): 8-45.

   \bibitem{BrennerJiangXu2009}
    M. Brenner, Y. Jiang, and Y. Xu, Multiparameter regularization for Volterra kernel identification via multiscale collocation methods, Advances in Computational Mathematics, 31 (2009): 421-455.
    
    
   \bibitem{Candes-Romberg-Tao:CPAM:06}
    E. J. Candès, J. Romberg, and T. Tao, Robust uncertainty principles: Exact signal reconstruction from highly  incomplete frequency information, IEEE Transactions on Information Theory, 52 (2006): 489-509.

  
    \bibitem{Candes2008}
    E. J. Candès and M. B. Wakin, An introduction to compressive sampling,
    IEEE Signal Process. Mag., 25 (2) (2008): 21-30

    \bibitem{ChenLuXuYang2008}
     Z. Chen, Y. Lu, Y. Xu, and H. Yang, Multi-parameter Tikhonov regularization for linear ill-posed operator equations, Journal of Computational Mathematics, 26 (1) (2008): 37-55.

	\bibitem{Chen2015}
	Z. Chen, C. Micchelli, and Y. Xu, Multiscale Methods for Fredholm Integral Equations, Cambridge University Press, Cambridge,  2015.

   
   \bibitem{Chong2017}
    E. Chong,  C. Han, and F.C. Park, Deep learning networks for stock market analysis and prediction: Methodology, data representations, and case studies, Expert Syst. Appl., 83 (2017): 187–205.	

	\bibitem{Cybenko1989}
	G. Cybenko, Approximation by superpositions of a sigmoidal function, Math. Control Signals Systems 2 (1989): 303-314.
	
	\bibitem{Cyr2020}
	E. C. Cyr, M. A. Gulian, R. G. Patel, M. Perego, and N. A. Trask, Robust training and initialization of deep neural networks: An adaptive basis viewpoint, vol. 107 of Proceedings of Machine Learning Research, Princeton University, Princeton, NJ, USA, 20–24 Jul 2020, PMLR, pp. 512–536, http://proceedings.mlr.press/v107/cyr20a.html
	
	\bibitem{Dahl2012}
	G. E. Dahl, D. Yu, L. Deng, and A. Acero, Context-dependent pre-trained deep neural networks
    for large-vocabulary speech recognition, IEEE Trans. Audio Speech Language Process., 20 (2012): 30-42.

   \bibitem{Daubechies1992}
   I. Daubechies, Ten Lectures on Wavelets, Society for Industrial and Applied Mathematics, Philadelphia, 1992.
   
	\bibitem{Daubechies2019}
   I. Daubechies, R. DeVore, S. Foucart, et al, Nonlinear Approximation and (Deep)  ReLU Networks, Constr. Approx (2021). https://doi.org/10.1007/s00365-021-09548-z


    \bibitem{Davis2020}
    D. Davis, D. Drusvyatskiy, S. Kakade, and J. D. Lee,  Stochastic subgradient method converges on tame functions, Found. Comput. Math., 20 (2020): 119–154.
    
    \bibitem{Bert2018}
    J. Devlin, M. W. Chang, K. Lee, and K. Toutanova, BERT: Pre-Training of deep bidirectional transformers for language understanding, 2018, https://arxiv.org/abs/1810.04805v1
    
    \bibitem{DeVore2021}
    R. DeVore, B. Hanin, and G. Petrova, Neural network approximation, Acta Numerica, 30 (2021): 327-444.
	
	\bibitem{Diefenthaler2021}
	M. Diefenthaler, A. Farhat, A. Verbytskyi, and Y. Xu, Deeply learning deep inelastic scattering kinematics, 2021, https://arxiv.org/abs/2108.11638v1
	
   
    \bibitem{Donoho:IEEEIT:06} D. Donoho, Compressive sensing, IEEE Transanctions on Information Theory, 52 (2006): 1289-1306.
   
   \bibitem{Driscoll:2014}
   T. A. Driscoll, N. Hale, and L. N. Trefethen, Chebfun Guide, 2014.

    \bibitem{Duchi2011}
    J. Duchi, E. Hazan, and Y. Singer, Adaptive subgradient methods for online learning and stochastic optimization,  Journal of Machine Learning Research, 12 (2011): 2121-2159.
    
	\bibitem{Friston2008}
    K. Friston, Hierarchical models in the brain, PLoS Comput. Biol, 4 (11) (2008): e1000211.
	
    \bibitem{Goodfellow2016}
    I. Goodfellow, Y. Bengio, and A. Courville, Deep Learning, MIT Press, 2016.
    
	\bibitem{Han_jentzen_E_2018}
	J. Han, A. Jentzen, and W. E, Solving high-dimensional partial differential equations using deep learning, Proceedings of the National Academy of Sciences 115 (2018): 8505-8510.
  
  \bibitem{He_Li_Xu_2020}
  J. He, L. Li, J. Xu, and C. Zheng, ReLU deep neural networks and linear finite elements, Journal of Computational Mathematics, 38 (2020): 502-527.
  
  \bibitem{Helton_2003}
   J. C. Helton and F.J. Davis, Latin hypercube sampling and the propagation of uncertainty in analyses of complex systems,
   Reliability Engineering \& System Safety, 81 (1) (2003): 23-69. 
   
   \bibitem{Hoefler_2021SIAMNews} T. Hoefler, D. A. Alistarh, N. Dryden, and T. Ben-Nun, The future of deep learning will be sparse, SIAM News, May 03, 2021.


   \bibitem{Hoefler_2021}
   T. Hoefler, D. A. Alistarh, T. Ben-Nun, N. Dryden,  and E. A. Peste,  Sparsity in deep learning: Pruning and growth for efficient inference and training in neural networks, Journal of Machine Learning Research, 22 (241) (2021): 1-124.
  
   \bibitem{Jung2020}
  J. Jung, K. Yoon, and P. Lee, Deep learned finite elements, Computer Methods in Applied Mechanics and Engineering, 372 (2020), 113401.

  \bibitem{Kingma2015}
  D. Kingma and J. Ba, Adam: A method for stochastic optimization, Proceedings of the 3rd International Conference on Learning Representations (ICLR 2015).

	\bibitem{Krizhevsky_2012}    
    A. Krizhevsky, I. Sutskever, and G. E. Hinton, ImageNet classification with deep convolutional
    neural networks, in Advances in Neural Information Processing Systems, Curran Associates, Inc.,
    2012, 1097-1105.	
	
	\bibitem{Lagaris1998}
	I. E. Lagaris, A. Likas, and D. I. Fotiadis, Artificial neural networks for solving ordinary and partial differential equations, IEEE Transactions on Neural Networks 9 (1998): 987-1000.
	
	\bibitem{Lagaris2000}
	I. E. Lagaris, A. C. Likas, and D. G. Papageorgiou, Neural-network methods for boundary value problems with irregular boundaries, IEEE Transactions on Neural Networks, 11 (2000): 1041-1049.

   \bibitem{Liu1995}
    T. Liu and K. Zumbrun, Nonlinear stability of an undercompressive shock for complex Burgers equation, Communications in Mathematical Physics, 168 (1) (1995): 163-186.

   \bibitem{LuShenXu2007}
    Y. Lu, L. Shen, and Y. Xu, Multi-Parameter regularization methods for high-resolution image reconstruction with displacement errors, IEEE Transactions on Circuits and Systems I, 54 (8) (2007): 1788-1799.
 
 \bibitem{Micchelli_Xu1994}
  C. A. Micchelli and Y. Xu, Using the matrix refinement equation for the construction of wavelets on invariant sets, Applied and Computational Harmonic Analysis, 1 (4) (1994): 391-401.
  
    \bibitem{Raissi2018}
    M. Raissi, Deep hidden physics models: deep learning of nonlinear partial differential equations, J. Mach. Learn. Res., 19 (1) (2018): 932-955.
	
	\bibitem{PINN_2019}
	M. Raissi, P. Perdikaris, and  G.E. Karniadakis,
	Physics-informed neural networks: A deep learning framework for solving forward and inverse problems involving nonlinear partial differential equations, Journal of Computational Physics, 378 (2019): 686-707.
	
	\bibitem{Serkin2000}
    V. N. Serkin and A. Hasegawa, Novel soliton solutions of the nonlinear Schrödinger equation model, Physical Review Letters, 85 (21) (2000): 4502.

    \bibitem{Unser2019}
    M. Unser, A representer theorem for deep neural networks,  Journal of Machine Learning Research, 20 (2019): 1-28.	   

   \bibitem{Xu_L0_2021}
   Y. Xu, Sparse regularization with the $\ell_0$ norm, 2021, arXiv:2111.08244. 
   
\bibitem{Xu_Ye2019}
   Y. Xu and Q. Ye, Generalized Mercer kernels and reproducing kernel Banach spaces, Memoirs of the
American Mathematical Society, 258 (2019): 1-122.
   
	\bibitem{XuZhang2021}
	Y. Xu and H. Zhang, Convergence of deep ReLU networks, 2021, arXiv:2107.12530.

	\bibitem{XuZhang_cnn2021}
	Y. Xu and H. Zhang, Convergence of deep convolutional neural networks, 2021, 	arXiv:2109.13542.
	
	\bibitem{Zhang-Xu-Zhang} H. Zhang, Y. Xu, and J. Zhang,
Reproducing kernel Banach spaces for machine learning, Journal of Machine Learning Research, 10 (2009): 2741-2775.
\end{thebibliography}
\end{document}